# Topological Approaches for Rough Continuous Functions with Suggested Applications

A. S. Salama
Faculty of Science, Tanta University, Tanta, Egypt
asalama@science.tanta.edu.eg

**Abstract:** In this paper, we purposed further study on rough functions and introduced some concepts based on it. We introduced and investigated the concepts of topological lower and upper approximations of near open sets and studied their basic properties. We defined and studied new topological neighborhood approach of rough functions. We generalized rough functions to topological rough functions by different topological structures. In addition, topological approximations of a function as a relation are defined and studied. Finally, we apply our approach of rough functions in finding the images of patient classification data using rough continuous functions.

***Keywords***: Topological Spaces; Rough sets; Interior Operator; Closure Operator; Continuous Functions; Ordered Relations.

**MSC: 54A05; 54B05; 54D35; 03B70**

## 1. Introduction

For a long time, general topologists faced many questions about the importance of abstract topological spaces. These questions were directed to them either from themselves or from others. The answers were always about the importance of general spaces in other branches of mathematics such as algebra, analysis, etc. Rough set theory is a recent approach for reasoning about data [1, 2, 3, 4, 5, 6, 15]. This theory depends on a certain topological structure that achieved a great success in many fields of real life applications [7, 8, 9, 10, 11, 12, 13]. Now the general topologists can give a more sensible answer to the above questions and can say, "Rough sets theory is a topological bridge from real life problems to computer science" [14, 16] .

Rough set theory was proposed as a new approach to processing of incomplete data. One of the aims of the rough set theory is a description of imprecise concepts. Suppose we are given a finite non-empty set $U$ of objects, called universe. Each object of $U$ is characterized by a description, for example a set of attributes values. In standard rough sets introduced by Pawlak an equivalence relation (reflexive, symmetric and transitive relation) on the universe of objects is defined based on their attribute values. In particular, this equivalence relation is constructed based on the equality relation on attribute values. Many real-world applications have both nominal and continuous attributes [17, 18, 19]. It was early recognized that standard rough set model based on the indiscernibility relation is well suited in the case of nominal attributes.

Many attempts were made to resolve limitations of this approach and many authors proposed interesting extension of the initial model (for example [20, 21, 22, 23, 24]). It was observed that considering a similarity relation instead of an indiscernibility relation is quite relevant. A binary relation forming classes of objects, which are identical or at least not noticeably different in terms of the available description, can represent the similarities between objects [25, 26, 27, 28, 29]. More recent approaches of rough set with its applications can found in [31,32,33]. Other applications of rough set theory in computer science (field of information retrievals) using topological generalizations can be found in [34.35.36.37].

In this paper, we purposed further study on rough functions and introduced some concepts based on rough functions. In Section 2, we give in more details the fundamentals of near open sets. The aim of Section 3 is to introduce and investigate the concepts of topological lower and upper approximations of near open sets and study their basic properties. The main goal of Section 4 is to spotlight on rough numbers. We aim in Section 5 to define and study new topological neighborhood approach of rough functions. Section 6 is devoted to generalize the concept of rough function to topological rough function by using different topological structures. Topological approximations of a function as a relation are defined and studied in Section 7. In Section 8 we suggested some applications of rough functions to information systems and give some applications of them in data retrieval. Finally, we concluded our work in Section 9.

## 2. Basic Concepts of Topological Near Open Sets

In this section, we give in more details the fundamentals of near open sets. The notions semi-open sets and semi-pre-open sets are introduced respectively in [30]. The generality of closed sets to generalized closed sets and to semi-generalized closed sets via semi-open sets are introduced in [30]. The complement of a semi-open (resp. semi-generalized closed) set is called a semi-closed (resp. semi-generalized open) set.

A topological space is a pair $(U,\tau)$ consisting of a set $U$ and family $\tau$ of subset of $U$ satisfying the following conditions:

(1) $\phi, U \in \tau$,
(2) $\tau$ is closed under arbitrary union.
(3) $\tau$ is closed under finite intersection.

The pair $(U,\tau)$ is called a topological space, the elements of $U$ are called points of the space, the subsets of $U$ belonging to $\tau$ are called open sets in the space, and the complement of the subsets of $U$ belonging to $\tau$ are called closed sets in the space; the family $\tau$ of open subsets of $U$ is also called a topology for $U$.

For a subset $A$ of a space $(U,\tau)$, $\overline{A}, A^o$ and $A^c$ denote the closure of $A$, the interior of $A$ and the complement of $A$ in $X$ respectively.

Let us recall the following definitions, which are useful in the sequel.

A subset $A$ of a topological space $(U,\tau)$ is called:

(1) Semi-open set if $A \subseteq \overline{(A^o)}$ and it is called a semi-closed set if $(\overline{A})^o \subseteq A$.

(2) Pre-open set if $A \subseteq (\overline{A})^o$ and it is called a pre-closed set if $\overline{(A^o)} \subseteq A$.

(3) $\alpha$ – open set if $A \subseteq ((\overline{A^o}))^o$ and it is called a $\alpha$ – closed set if $\overline{((\overline{A})^o)} \subseteq A$.

(4) semi-pre-open set ($\beta$ – open) if $A \subseteq \overline{((\overline{A})^o)}$ and it is called a semi-pre-closed set ($\beta$ – closed) if $((\overline{A^o}))^o \subseteq A$.

(5) Regular-open set if $A \subseteq (\overline{A})^o$ and it is called a regular-closed set if $\overline{(A^o)} = A$.

(6) Semi-regular set if it both semi-open and semi-closed in $(U,\tau)$.

(7) $\delta$ – closed set if $A = \overline{\delta(A)}$, where $\overline{\delta(A)} = \{x \in U : (\overline{G})^o \cap A \neq \varphi, x \in G, G \in \tau\}$.

The semi-closure (resp. $\alpha$-closure, semi-pre-closure) of a subset $A$ of $(U,\tau)$ is the intersection of all semi-closed (resp. $\alpha$–closed, semi-pre-closed) sets that contains $A$ and is denoted by $S(\overline{A})$ (resp. $\alpha(\overline{A}), sp(\overline{A})$). The union of all semi-open subsets of $U$ is called the semi-interior of $A$ and is denoted by $s(A^o)$.

A subset $A$ of a topological space $(U,\tau)$ is called:

(1) Generalized closed set if $\overline{A} \subseteq G$ whenever $A \subseteq G$ and $G \in \tau$.

(2) Semi-generalized closed set if $s(\overline{A}) \subseteq G$ whenever $A \subseteq G$ and $G$ is semi-open set in $(U,\tau)$.

The complement of a $sg$ – closed set is called a $sg$ – open set.

(3) Generalized semi-closed set if $s(\overline{A}) \subseteq U$ whenever $A \subseteq G$ and $G \in \tau$.

(4) $\alpha$ – Generalized closed set if $\alpha(\overline{A}) \subseteq G$ whenever $A \subseteq G$ and $G \in \tau$.

(5) Generalized $\alpha$ – closed set if $\alpha(\overline{A}) \subseteq G$ whenever $A \subseteq G$ and $G$ is $\alpha$ – open in $(U, \tau)$.

(6) $g\alpha^{**}$ -closed set if $\overline{A} \subseteq \overline{(G^o)}$ whenever $A \subseteq G$ and $G$ is $\alpha$ – open in $(U, \tau)$.

**3. Topological Near Open Approach of Rough Approximations**

In this section, we introduce and investigate the concepts of topological lower and upper approximations of near open sets and study their basic properties.

Let $(U, \tau)$ be a topological space. If $X \subseteq U$, then:

(1) Semi-lower approximation of $X \subseteq U$, $\underline{X}_s = \bigcup\{G : G \in Semi(U), G \subseteq X\}$ where $Semi(U)$ is the set of all semi-open sets in $(U, \tau)$. Semi-upper approximation of $X \subseteq U$, $\overline{X}_s = \bigcap\{F : F \in CSemi(U), F \cap X \neq \varphi\}$ where $CSemi(U)$ is the set of all semi-closed sets in $(U, \tau)$.

(2) Pre-lower approximation of $X \subseteq U$, $\underline{X}_P = \bigcup\{G : G \in \Pr e(U), G \subseteq X\}$ where $\Pr e(U)$ is the set of all pre-open sets in $(U, \tau)$. Pre-upper approximation of $X \subseteq U$, $\overline{X}_P = \bigcap\{F : F \in C \Pr e(U), F \cap X \neq \varphi\}$ where $C \Pr e(U)$ is the set of all pre-closed sets in $(U, \tau)$.

(3) $\alpha$-lower approximation of $X \subseteq U$, $\underline{X}_\alpha = \bigcup\{G : G \in \alpha(U), G \subseteq X\}$, where $\alpha(U)$ is the set of all $\alpha$- open sets in $(U, \tau)$. $\alpha$-upper approximation of $X \subseteq U$, $\overline{X}_\alpha = \bigcap\{F : F \in C\alpha(U), F \cap X \neq \varphi\}$, where $C\alpha(U)$ is the set of all $\alpha$-closed sets in $(U, \tau)$.

(4) $\beta$-lower approximation of $X \subseteq U$, $\underline{X}_\beta = \bigcup\{G : G \in \beta(U), G \subseteq X\}$ where $\beta(U)$ is the set of all $\beta$-open sets in $(U, \tau)$. $\beta$-upper approximation of $X \subseteq U$, $\overline{X}_\beta = \bigcap\{F : F \in C\beta(U), F \cap X \neq \varphi\}$ where $C\beta(U)$ is the set of all $\beta$-closed sets in $(U, \tau)$.

(5) Regular-lower approximation of $X \subseteq U$, $\underline{X}_{Reg} = \bigcup\{G : G \in \Re e g(U), G \subseteq X\}$ where $\Re e g(U)$ is the set of all regular-open sets in $(U, \tau)$. Regular-upper approximation of $X \subseteq U$, $\overline{X}_{Reg} = \bigcap\{F : F \in C \Re e g(U), F \cap X \neq \varphi\}$ where $C \Re e g(U)$ is the set of all regular-closed sets in $(U, \tau)$.

(6) Semi-regular lower approximation of $X \subseteq U$, $\underline{X}_{SReg} = \bigcup\{G : G \in S \Re e g(U), G \subseteq X\}$, where $S \Re e g(U)$ is semi-closed sets in $(U, \tau)$ respectively. Semi-regular upper approximation of $X \subseteq U$, $\overline{X}_{SReg} = \bigcap\{F : F \in CS \Re e g(U), F \cap X \neq \varphi\}$ where $CS \Re e g(U)$ is the set of all semi-closed sets in $(U, \tau)$.

Motivation for topological rough set theory has come from the need to represent subsets of a universe in terms of topological classes of the topological base generated by the general binary relation defined on the universe. That base characterizes a topological space, called topological approximation space $App_\tau = (U, R, \tau_R)$. The topological classes of $R$ are also known as the topological granules, topological elementary sets or topological blocks; we will use $G_{xm} \in \tau$, to denote the topological class containing $x \in U$. In the topological approximation space, we consider two operators $\underline{R}_m(X) = \{x \in U : G_{xm} \subseteq X\}$ and $\overline{R}_m(X) = \{x \in U : G_{xm} \cap X \neq \varphi\}$ called the topological lower approximation and topological upper approximation of $X \subseteq U$ respectively. Also let $POS_m(X) = \underline{R}_m(X)$ denote the topological positive region of $X \subseteq U$, $NEG_m(X) = U - \overline{R}_m(X)$ denotes the topological negative region of $X \subseteq U$ and $BON_m(X) = \overline{R}_m(X) - \underline{R}_m(X)$ denote the topological borderline region of $X \subseteq U$.

The degree of topological completeness can characterize by the topological accuracy measure, in which $|X|$ represents the cardinality of set $X \subseteq U$ as follows:

$$\alpha_m(X) = |\underline{R}_m(X)|/|\overline{R}_m(X)|, X \neq \varphi.$$

We define here the semi-rough pairs as an example of topological rough sets and we studying their properties. You can use any type of the above near open sets as another example.

The semi topological class on a topological approximation space $App_\tau = (U, R, \tau_R)$ is determined by $(\underline{X}_s, \overline{X}_s) = \{A \subset U : \underline{X}_s \subset A \subset \overline{X}_s\}$. A subset $X \subset U$ is said to be semi-dense (semi-co-dense) if $\overline{X}_s = U$ ($\underline{X}_s = \varphi$). By semi-rough pair on $App_\tau = (U, R, \tau_R)$ we mean any pair $(P, Q)$ where $P, Q \subseteq U$ satisfying the conditions:

(Semi-1) $P$ is semi-open set in $\tau_R$.
(Semi-2) $Q$ is semi-closed set in $\tau_R$.
(Semi-3) $P \subset Q$.
(Semi-4) there is a subset $S \subset U$ such that:
　1- $S_s^o = \varphi$,
　2- $S \subset Q - \overline{P}_s$,
　3- $Q - \overline{P}_s \subset \overline{S}_s$.

**Lemma 3.1** For any subset $A \subset U$ in the topological approximation space $App_\tau = (U, R, \tau_R)$, the pair $(\underline{A}_s, \overline{A}_s)$ is a semi-rough pair on $App_\tau = (U, R, \tau_R)$ in which every semi-open set in $U$ is a semi closed set.

**Proof:** Let $P = A_s^o$ and $Q = \overline{A}_s$. Then the conditions from (Semi-1) to (Semi-3) are directly satisfied. Now we need to prove condition (Semi-4). Define $S = A - \overline{P}_s$ then we have:

1- If $O \subseteq S$ is a semi-open set, hence $O \subseteq A$ that gives $O \cap P = \varphi$ which is a contradiction, hence $O$ is not contained in $A$. Then must $S = \varphi$ which give $S_s^o = \varphi$.

2- Since $S = A - \overline{P}_s$, $A \subseteq \overline{A}_s$, then $A \subset Q$. Then we have $S \subset Q - \overline{P}_s$.

3- Let $x \in Q - \overline{P}_s$, $Q = \overline{A}_s$, this means that $x \in A$ or $x \notin A$. If $x \in A$, then for every semi-open set $O$ and $x \in O$ such that $O \cap A \neq \varphi$ implies that $O \cap S \neq \varphi$ and we have $x \in \overline{S}_s$, then $Q - \overline{P}_s \subset \overline{S}_s$. If $x \notin A$, then there is a semi-open set $O'$ and $x \in O'$. Now $O' - \overline{P}_s = O' \cap [\overline{P}_s]^c$ is a semi-open set contains $x$, and $x \in \overline{A}_s$, then there exist a point $y \in A$ such that $y \in O' - \overline{P}_s$, hence $y \in O' \cap S$, therefore $O' \cap S \neq \varphi$, hence $x \in \overline{S}_s$. Then we have the result $Q - \overline{P}_s \subset \overline{S}_s$. □

**Lemma 3.2** For any semi-rough pair $(P, Q)$ in $App_\tau = (U, R, \tau_R)$ in which every semi-open subset is semi-closed, there are subsets $A, B \subseteq U$ such that $P = \overline{A}_s$ and $Q = \overline{B}_s$.

**Proof:** Let $(P, Q)$ be a semi-rough pair and let $S$ be any subset satisfying condition (Semi-4). Define $A = P \cup S$, then $P \subset A$ hence $P \subset A_s^o$. If $O \subset A$ is a semi-open set then $O - \overline{P}_s = O \cap [\overline{P}_s]^c$ is another semi-open set contained in $A$. Since $O \subset A = P \cup S$, $P \subset \overline{P}_s$, then $O \subset \overline{P}_s \cup S$ and we have $O - \overline{P}_s \subset [(\overline{P}_s \cup S) - \overline{P}_s] = S$. Therefore $O - \overline{P}_s$ is a semi-open set contained in $S$ which means $O \subset \overline{P}_s$. Since $O \subset P \cup S$, it follows that $O \subset P$ and this prove that $P = \overline{A}_s$.

Now we have $\overline{P}_s \cup \overline{S}_s \supset \overline{P}_s \cup [Q - \overline{P}_s] = Q$. Also, $B = P \cup S \supset \overline{P}_s \cup (Q - \overline{P}_s) = Q$ and hence $\overline{B}_s \subset \overline{Q}_s = Q$. Then we have $Q = \overline{B}_s$. □

**Theorem 3.1** For any topological subspace $(X, \tau^*), X \subseteq U$ of the topological approximation space $App_\tau = (U, R, \tau_R)$ the function $f : (X, \tau^*) \to (U, \tau_R)$ that defined by $f(A) = (A_s^o, \overline{A}_s), A \in \tau^*$ is bijection.

**Proof:** First, we will prove that the function is onto: For any semi-rough pair $(A_s^o, \overline{A}_s)$ in $App_\tau = (U, R, \tau_R)$ then there exist $A \in \tau^*$ such that $f(A) = (A_s^o, \overline{A}_s)$. Second, for the proof that $f$ is one to one: if $f(A_1) = f(A_2)$ then $(A_{1s}^o, \overline{A}_{1s}) = (A_{2s}^o, \overline{A}_{2s})$ which imply to $A_{1s}^o = A_{2s}^o$ and $\overline{A}_{1s} = \overline{A}_{2s}$ and $A_1 \approx A_2$. □

**4. Approach of Rough Numbers**

Rough set theory depending essentially on equivalence relation $R$ that have a class of elementary blocks denoted by $U/R$. The relation $R$ on $U$ have three different cases, the first is the equivalence case and in this case $U/R$ partitioned the universal set. The second case when $R$ is non equivalence but their elementary blocks cover $U$, then $U/R$ denoted by $C(U)$ in some works. Finally, third case that $R$ is any binary relation that neither partition nor covering of $U$, then some researchers use $U/R$ as a subbase of a topology on $U$ and denoted $U/R$ in this case by $\beta(U)$. We used $R(U)$ in this section to consider any one of the above cases and we added more examples and explanations of $R(U)$ in the section.

Let $U$ be the set of non-negative real numbers, consider the numbers $x_1, x_2, \ldots, x_i \in U$ and define the order relation $R = \{x_n R x_m : x_n < x_m, n < m < i, \forall n, m, i \in U\}$. The relation $R$ defines the partition $R(U)$ of $U$ by $R(U) = \{x_i, (x_i, x_{i+1}), x_{i+1}, \forall i \in U\}$, where $(x_i, x_{i+1})$ denote open intervals in $U$. The relation $R$ is called a classification of $U$ and the ordered pair $Appr(U) = (U, R)$ is an approximation space, where $R$ is the binary relation associated with $R(U)$.

Let $Appr(U) = (U, R)$ be an approximation space. By $R(x)$ in $Appr(U) = (U, R)$ we denote the block of the partition $R(U)$ containing $x$, in particular if $x \in R(x)$, we have $R(x) = \{x\}$. The approximation space $Appr(U) = (U, R)$ is selective if for all $x \in U, [x]_R = \{x\} = R(x)$

Let $Appr(U) = (U, R)$ be an approximation space, by $\overline{R}(x)$ we denote the closed interval $[0, x]$ for $x \in U$.

For any $x \in U$, the R-lower and the R-upper approximations of $\overline{R}(x)$ in the approximation space $Appr(U) = (U, R)$ are defined respectively by:

$$\underline{Appr(U)[\overline{R}(x)]} = \{y \in U : R(y) \subseteq \overline{R}(x)\}$$
$$\overline{Appr(U)[\overline{R}(x)]} = \{y \in U : R(y) \cap \overline{R}(x) \neq \varphi\}$$

The approximations of the closed interval $\overline{R}(x) = [0, x]$ can be understood as the approximations of the real number $x$. The number $x$ is a rough number if $\underline{Appr(U)[\overline{R}(x)]} \neq \overline{Appr(U)[\overline{R}(x)]}$, otherwise it is exact number.

**Example 4.1** Let $U$ be the set of non-negative real numbers. We define the infinite partition $R(U) = \{0, (0,1), 1, (1,2), 2, \ldots, n, (n, n+1), n+1, \ldots\}$ and hence, $Appr(U) = (U, R)$ is an approximation space. Also, for any integer number $x \in U$, we have $R(x) = \{x\}$ and for non-integers, $R(x) = (x_i, x_{i+1})$ and $x \in (x_i, x_{i+1})$, $x_i, x_{i+1} \in U$. Then every non integer number is a rough number in $Appr(U) = (U, R)$.

**Example 4.2** Let $R^1$ be the set of real numbers and let $(c,d) \subseteq R^1$ be an open interval. We convert the interval $(c,d)$ by the sequence $S = \{c = x_0, x_1, ..., x_n = d\}$ of real numbers such that $c = x_0 < x_1 < ... < x_{n-1} < x_n = d$. we assume that $0 \in S$. The partition generated by the sequence $S$ is given by: $R(S) = \{\{c\}, (c, x_1), \{x_1\}, (x_1, x_2), ..., \{x_{n-1}\}, (x_{n-1}, x_n), \{x_n\}\}$. Then the pair $(R^1, R(S))$ is an approximation space. We say that the number $x \in R^1$ is rough in $(R^1, R(S))$ if $\underline{Appr}(R^1)[\overline{R}(x)] \neq \overline{Appr}(R^1)[\overline{R}(x)]$.

## 5. Topological Neighborhood Approach of Rough Continuity

Let $X$ and $Y$ be two subsets of a universe $U$, and let $Appr(X) = (X, S)$ and $Appr(Y) = (Y, P)$ be two approximation spaces, where $S$ and $P$ are binary relations on $X$ and $Y$ respectively. We define two subsets $S_r(x) = \{y \in X : (x, y) \in S\}$ and $S_l(x) = \{y \in X : (y, x) \in S\}$ of $X$ (also two subsets $P_r(x) = \{y \in Y : (x, y) \in P\}$ and $P_r(x) = \{y \in Y : (y, x) \in P\}$ of $Y$) are called right and left neighborhoods of an element $x \in X$. We define now two topologies on $X$ and on $Y$, respectively, using the intersection of the right and left neighborhoods $S_{r \cap l}(x) = S_r(x) \cap S_l(x)$ and $P_{r \cap l}(x) = P_r(x) \cap P_l(x)$ as follows:
$$\tau_X = \{A \subseteq X : \forall a \in A, S_{r \cap l}(a) \subseteq A\},$$
$$\tau_Y = \{B \subseteq Y : \forall b \in B, P_{r \cap l}(b) \subseteq B\}.$$

The rough approximations using these topologies are defined as follows:
$$\underline{P}_{\tau_Y}(B) = \bigcup\{G' \in \tau_Y : G' \subseteq B\},$$
$$\underline{S}_{\tau_X}(A) = \bigcup\{G \in \tau_X : G \subseteq A\},$$
$$\overline{S}_{\tau_X}(A) = \bigcap\{F : F^c \in \tau_X : A \subseteq F\},$$
$$\overline{P}_{\tau_Y}(B) = \bigcap\{F' : F'^c \in \tau_Y : B \subseteq F'\}.$$

The function $f:(X, \tau_X) \longrightarrow (Y, \tau_Y)$ is called a rough function on $X$ if the image of each rough set in $X$ is rough in $Y$.

Namely, the function $f$ is totally rough iff all subsets $A \subset X, A \neq \varphi$, such that $\underline{S}_{\tau_X}(A) \neq \overline{S}_{\tau_X}(A)$ then $\underline{P}_{\tau_Y}(f(\underline{S}_{\tau_X}(A))) \neq \overline{P}_{\tau_Y}(f(\overline{S}_{\tau_X}(A)))$ in $Y$.

The function $f$ is possibly rough iff some subsets $A \subset X, A \neq \varphi$, such that $\underline{S}_{\tau_X}(A) \neq \overline{S}_{\tau_X}(A)$ then $\underline{P}_{\tau_Y}(f(\underline{S}_{\tau_X}(A))) \neq \overline{P}_{\tau_Y}(f(\overline{S}_{\tau_X}(A)))$ in $Y$.

Finally, the function $f$ is exact iff all subsets $A \subset X, A \neq \varphi$, such that $\underline{S}_{\tau_X}(A) = \overline{S}_{\tau_X}(A)$ then $\underline{P}_{\tau_Y}(f(\underline{S}_{\tau_X}(A))) = \overline{P}_{\tau_Y}(f(\overline{S}_{\tau_X}(A)))$ in $Y$.

The function $f:(X, \tau_X) \longrightarrow (Y, \tau_Y)$ is a topological rough, continuous function on $X$ as the following:

1) The function $f$ is topological totally rough, continuous if for all subsets $A \subset Y, A \neq \varphi$, if $(A)^o_{\tau_Y} \subseteq \overline{(A)}_{\tau_Y}$ then $(f^{-1}(\overline{(A)}_{\tau_Y}))^o_{\tau_X} \subseteq \overline{(f^{-1}((A)^o_{\tau_Y}))}_{\tau_X}$ in $X$.

2) The function $f$ is topological possibly rough, continuous if for some subsets $A \subset Y, A \neq \varphi$, if $(A)^o_{\tau_Y} \subseteq \overline{(A)}_{\tau_Y}$ then $(f^{-1}(\overline{(A)}_{\tau_Y}))^o_{\tau_X} \subseteq \overline{(f^{-1}((A)^o_{\tau_Y}))}_{\tau_X}$ in $X$.

3) Finally, the function $f$ is topological exact continuous if for all subsets $A \subset Y, A \neq \varphi$, if $(A)^o_{\tau_Y} = \overline{(A)}_{\tau_Y}$ then $(f^{-1}(\overline{(A)}_{\tau_Y}))^o_{\tau_X} = \overline{(f^{-1}((A)^o_{\tau_Y}))}_{\tau_X}$ in $X$.

**Example 5.1** Let $(X, \tau_X)$ and $(Y, \tau_Y)$ be topological spaces, where $X = \{a, b, c\}$ and $\tau = \{X, \varphi, \{a\}, \{b\}, \{a,b\}\}$ and $Y = \{1, 2, 3\}$, $\tau_Y = \{Y, \varphi, \{1\}, \{2\}, \{1,2\}\}$. Let $f: X \longrightarrow Y$ be a map defined by $f(a) = 1$, $f(b) = 2$ and $f(c) = 3$, then our results are given in Table 1 below.

| Subsets of $Y$ / Our measures | $\{1\}$ | $\{2\}$ | $\{3\}$ | $\{1,2\}$ | $\{1,3\}$ | $\{2,3\}$ | $Y$ |
|---|---|---|---|---|---|---|---|
| $(A)^o_{\tau_Y}$ | $\{1\}$ | $\{2\}$ | $\varphi$ | $\{1,2\}$ | $\{1\}$ | $\{2\}$ | $Y$ |
| $\overline{(A)}_{\tau_Y}$ | $\{1,3\}$ | $\{2,3\}$ | $\{3\}$ | $Y$ | $\{1,3\}$ | $\{2,3\}$ | $Y$ |
| $f^{-1}((A)^o_{\tau_Y})$ | $\{a\}$ | $\{b\}$ | $\varphi$ | $\{a,b\}$ | $\{a\}$ | $\{b\}$ | $X$ |
| $f^{-1}(\overline{(A)}_{\tau_Y})$ | $\{a,c\}$ | $\{b,c\}$ | $\{c\}$ | $X$ | $\{a,c\}$ | $\{b,c\}$ | $X$ |
| $(f^{-1}(\overline{(A)}_{\tau_Y}))^o_{\tau_Y}$ | $\{a\}$ | $\{b\}$ | $\varphi$ | $X$ | $\{a\}$ | $\{b\}$ | $X$ |
| $\overline{(f^{-1}((A)^o_{\tau_Y}))}_{\tau_X}$ | $\{a,c\}$ | $\{b,c\}$ | $\varphi$ | $X$ | $\{a,c\}$ | $\{b,c\}$ | $X$ |

Table 1: Calculations of topological rough continuous functions

Then, according to Table 1 the function $f$ is topological totally rough continuous function.

**Proposition 5.1** Let $(X, \tau_X)$ and $(Y, \tau_Y)$ be topological spaces and let $f: X \to Y$ be a function. The following are equivalent:

1. $f$ is rough continuous.
2. For every $F \subseteq Y$, $f^{-1}(\overline{(F)}_{\tau_Y}) \in \tau^c_X$.
3. For every $x \in X$, $f$ is rough continuous at $x$.
4. For every $A \subseteq X$, $f(\overline{A}_{\tau_X}) \subseteq \overline{f(A)}_{\tau_Y}$.

**Proof:** We will use the sequence 3) implies 1) implies 4) implies 2) implies 3) to prove the equivalence of the proposition.

3) implies 1): suppose a non-empty open set $V \in \tau_Y$, for a fixed point $x \in f^{-1}(V)$ we have $f(x) \in V$. But since $f$ is rough continuous at $x$, then there exist an open set $G_x \subseteq X$ such that $f(G_x) \subset V$ and $(f(G_x))^o_{\tau_Y} \subseteq \overline{(f(G_x))}_{\tau_Y}$, then we have $(f^{-1}(\overline{(f(G_x))}_{\tau_Y}))^o_{\tau_X} \subseteq \overline{(f^{-1}((f(G_x))^o_{\tau_Y}))}_{\tau_X}$ and $G_x \in f^{-1}(V)$, this give that $f$ is rough continuous.

1) implies 4): suppose that $f$ is rough continuous and let $A \subseteq X$. Let $x \in \overline{A}_{\tau_X}$. Let an open set $V \in \tau_Y$ such that $x \in f^{-1}(V)$. Then by the definition of rough upper approximation $f^{-1}(V) \cap A \neq \varphi$. Let $x' \in f^{-1}(V) \cap A$ then $f(x') \in V \cap f(A)$. Then we have $V \cap f(A) \neq \varphi$. Then we have $f(\overline{A}_{\tau_X}) \subseteq \overline{f(A)}_{\tau_Y}$.

4) implies 2): fix a closed subset $F \subseteq Y$, let $A = f^{-1}(F)$, we will prove that $A = \overline{A}_{\tau_X}$. But each subset is contained in its upper approximation, $A \subset \overline{A}_{\tau_X}$. Now we will prove that $\overline{A}_{\tau_X} \subset A$. Let $x \in \overline{A}_{\tau_X}$ then using 4) we have $f(x) \in f(\overline{A}_{\tau_X}) \subseteq \overline{f(A)}_{\tau_Y} \subseteq \overline{F}_{\tau_Y} = F$, hence $f(x) \in F$ or $x \in f^{-1}(F) = A$. Then we have $f^{-1}(\overline{(F)}_{\tau_Y}) \in \tau_X^c$.

2) implies 3): let $x \in X$ and $V \in \tau_Y$ an open set containing $f(x)$. Then $Y - V$ is a closed set and $f^{-1}(Y - V)$ is a closed set in $X$ does not contain the point $x$. But $x \in X - (f^{-1}(Y - V))$. Then there exists an open set $G$ containing $x$ such that $x \in G \subseteq X - (f^{-1}(Y - V))$ then $f(G) \subseteq f(X - (f^{-1}(Y - V))) = f(X) - (Y - V) \subseteq V$. Then $f$ is rough continuous at $x$. □

**Theorem 5.1** Suppose that $(\tau_i)_X, i = 1, 2, 3, ...$ be a family of topologies defined on $X$. Let $f : X \to Y$ be a rough continuous function for every $\tau_i, \forall i$ where $(Y, \tau_Y)$ is a topological space. Then $f$ is rough continuous function with respect to the topology $\tau_X = (\bigcap_i \tau_i)_X$.

**Proof:** let $G \in \tau_Y$, then $(G)_{\tau_Y}^o \subseteq \overline{(G)}_{\tau_Y}$, since $f$ is rough continuous function for every $\tau_i, \forall i$ then $(f^{-1}(\overline{(G)}_{\tau_Y}))^o_{\tau_i} \subseteq \overline{(f^{-1}((G)_{\tau_Y}^o))}_{\tau_i}, \forall i$. Then we have $(f^{-1}(\overline{(G)}_{\tau_Y}))^o_{\tau_X} \subseteq \overline{(f^{-1}((G)_{\tau_Y}^o))}_{\tau_X}$ in $\tau_X = (\bigcap_i \tau_i)_X$, hence $f$ is rough continuous function with respect to the topology $\tau_X = (\bigcap_i \tau_i)_X$. □

**Theorem 5.2** Let $f_i : X \to (Y_i, \tau_i)$ be a family of functions. Suppose that $\tau_X$ is the topology on $X$ generated by the class $\beta = \bigcup_i \{f_i^{-1}(G) : G \in \tau_i\}$ then:

1- $f_i$ is rough continuous for each $i$.
2- If $\tau_X^*$ is the intersection of all topologies on $X$ such that $f_i$ is rough continuous for each $i$. Then $\tau_X = \tau_X^*$.
3- $\tau_X$ is the coarser topology on $X$ gives that $f_i$ is rough continuous for each $i$.
4- The class $\beta = \bigcup_i \{f_i^{-1}(G) : G \in \tau_i\}$ is a subbase of $\tau_X$.
5- The function $g : Y \to X$ is rough continuous if and only if $f_i \circ g$ is rough continuous.

**Proof: Part 1)** for each function $f_i : X \to (Y_i, \tau_i)$ if $F \in \tau_i$ then $(F)_{\tau_i}^o \subseteq \overline{(F)}_{\tau_i}$, and $f_i^{-1}(F) \in \beta$. But $\beta \subseteq \tau_i$ then $f_i^{-1}(F) \in \tau_i$, hence $(f_i^{-1}(\overline{(F)}_{\tau_i}))^o_{\tau_X} \subseteq \overline{(f_i^{-1}((F)_{\tau_i}^o))}_{\tau_X}$ then we have the result.

**Part 2)** We can easily prove that $\beta \subseteq \tau_X^*$, but the topology $\tau_X$ is generated by $\beta$ then $\tau_X \subseteq \tau_X^*$. Otherwise, $\tau_X$ is one of the topologies that make the functions $f_i$ are rough continuous. Then we have $\tau_X^* \subseteq \tau_X$, hence $\tau_X = \tau_X^*$.

**Part 3)** Is obvious by proof of Part 2).

**Part 4)** since any collection of subsets of $X$ is a subbase of a topology on $X$. Then $\beta$ is a subbase of the topology $\tau_X$.

Part 5) if the function $g:Y \to X$ is rough continuous then all functions $f_i \circ g$ are rough continuous. Otherwise, let $f_i \circ g$ are rough continuous and let $G \in \beta$, then there exists a subset $H \in \tau_i$ such that $G = f_i^{-1}(H)$. But $g^{-1}(G) = g^{-1}(f_i^{-1}(H)) = (f_i \circ g)^{-1}(H)$. Now we have $(H)_{\tau_i}^o \subseteq \overline{(H)}_{\tau_i}$, then $(f_i \circ g)^{-1}((\overline{H})_{\tau_i}))^o{}_{\tau_X} \subseteq \overline{((f_i \circ g)^{-1}((H)_{\tau_i}^o))}_{\tau_X}$. Then $f_i \circ g$ is rough continuous.

□

## 6. Minimal Neighborhood Approach for Rough Continuity

We propose to generalize the concept of rough function to topological rough function by using topological structures. The topological spaces with rough sets are very useful in the field of digital topology, which is applied widely in the image processing in computer sciences.

Let $(X,\tau)$ be a topological space and $x \in X$. Then we define:
$N_{\min}(x) = \bigcap\{N \subseteq X : x \in G \subseteq N, \forall G \in \tau, x \in G\}$ which is called the minimal neighborhood containing the point $x$ with respect to the topology $\tau$ on $X$. Let $(X,\tau)$ be a topological space, for any element $x \in X$, we define the subset $\overline{N}_{\min}(x)$ which is the closure of $N_{\min}(x)$ in $(X,\tau)$.

If $f:(X,\tau) \longrightarrow (Y,\tau^*)$ is a function between two topological spaces $(X,\tau)$ and $(Y,\tau^*)$, we define the functions $f_{\min}:(X,\tau) \longrightarrow (Y,\tau^*)$, by:
$f_{\min}(x) = \bigcap\{M \subseteq Y : f(x) \in G' \subseteq M, \forall G' \in \tau^*, x \in G'\}$ for every $x \in X$.

Let $f:(X,\tau) \longrightarrow (Y,\tau^*)$ be a function, where $X$ and $Y$ are topological spaces. The function $f$ is called a topological rough function on $X$ if and only if $(N_{\min}(x))_\tau^o \neq \overline{(N_{\min}(x))}_\tau$ for every $x \in X$. Also, $f$ is a topological rough function on $Y$ if $(f_{\min}(x))_{\tau^*}^o \neq \overline{(f_{\min}(x))}_{\tau^*}$ for every point $f(x)$ in $Y$.

**Example 6.1** Let $(X,\tau)$ and $(Y,\tau^*)$ be topological spaces, where $X = \{a,b,c\}$ and $\tau = \{X,\varphi,\{a\},\{a,b\}\}$ and $Y = \{1,2,3\}$, $\tau^* = \{Y,\varphi,\{1\},\{2\},\{1,2\}\}$. Let $f: X \longrightarrow Y$ be a map defined by $f(a) = 2$, $f(b) = 1$ and $f(c) = 3$, then
$N_{\min}(a) = \{a\}$, $N_{\min}(b) = \{a,b\}$, $N_{\min}(c) = X$, $f_{\min}(a) = \{2\}$, $f_{\min}(b) = \{1\}$, $f_{\min}(c) = Y$,
Then we have:
$$(N_{\min}(a))_\tau^o = \{a\}, \text{ and } \overline{(N_{\min}(a))}_\tau = X,$$
$$(N_{\min}(b))_\tau^o = \{a,b\}, \text{ and } \overline{(N_{\min}(b))}_\tau = X,$$
$$(N_{\min}(c))_\tau^o = X, \text{ and } \overline{(N_{\min}(c))}_\tau = X,$$
Also,
$$(f_{\min}(a))_{\tau^*}^o = \{2\}, \text{ and } \overline{(f_{\min}(a))}_{\tau^*} = \{2,3\},$$
$$(f_{\min}(b))_{\tau^*}^o = \{1\}, \text{ and } \overline{(f_{\min}(b))}_{\tau^*} = \{1,3\},$$
$$(f_{\min}(c))_{\tau^*}^o = \{3\}, \text{ and } \overline{(f_{\min}(c))}_{\tau^*} = \{3\},$$
Then the function $f$ is not topological rough function on $X$ and on $Y$.

A function $f:(X,\tau) \longrightarrow (Y,\tau^*)$ is said to be a topological roughly continuous at the point $x \in X$ if and only if $f^{-1}(N_{\min}(f(x))) \subseteq N_{\min}(x)$, and it is a topological roughly continuous on $X$ if it is a topological roughly continuous at every point $x \in X$.

**Example 6.2** Let $f:(X,\tau) \longrightarrow (Y,\tau^*)$ be a function defined by $f(a) = 2$, $f(b) = f(d) = 3$ and $f(c) = 4$, where $X = \{a,b,c,d\}$ and $Y = \{1,2,3,4\}$ with
$$\tau = \{X,\varphi,\{a\},\{a,b\},\{a,b,c\}\},$$

and
$$\tau^* = \{Y, \varphi, \{1\}, \{2\}, \{1,2\}, \{2,3,4\}\}.$$
Then $f$ is a topological rough function on $X$ and

$N_{min}(a) = \{a\}$, but $N_{min}(f(a)) = N_{min}(2) = \{2\}$ then, $f^{-1}(N_{min}(2)) = \{a\}$

$N_{min}(b) = \{a,b\}$, but, $N_{min}(f(b)) = N_{min}(3) = \{2,3,4\}$ then, $f^{-1}(N_{min}(3)) = X$

$N_{min}(c) = \{a,b,c\}$, but, $N_{min}(f(c)) = N_{min}(4) = \{2,3,4\}$ then, $f^{-1}(N_{min}(4)) = X$

$N_{min}(d) = X$, but, $N_{min}(f(d)) = N_{min}(3) = \{2,3,4\}$ then, $f^{-1}(N_{min}(3)) = X$

then $f^{-1}(N_{min}(f(x))) \subseteq N_{min}(x)$ for every $x \in X$,

then $f$ is a topological rough, continuous function on $X$.

### 7. Topological Approximations of a Function as a Relation

The function $f: X \to Y$ is a relation from $X$ to $Y$ when it satisfies the conditions:

(i) $Dom(f) = X$

(ii) If $(x,y) \in f$ and $(x,z) \in f$, then $y = z$.

If $X = Y$, we say $f$ is a function on $X$. By this way any function $f: X \to Y$ can completely represented by its graph $G(f) = \{(x, f(x)) : x \in X\}$.

Let $f:(U_1, R_1) \longrightarrow (U_2, R_2)$ be any function, where $A_1 = (U_1, R_1)$ and $A_2 = (U_2, R_2)$ are approximation spaces, such that $R_1$ and $R_2$ are any binary relations on $U_1$ and $U_2$ respectively. We define the relation $R = R_1 \times R_2$ such that $R(x) = R_1(x) \times R_2(x)$ is the blocks of $U_1 \times U_2$. For the function, $G(f) = \{(x, f(x)) : x \in U_1\}$ we define the approximations

$$\underline{R}(G(f)) = \cup \{G \subseteq R(x) : G \in G(f)\},$$
$$\overline{R}(G(f)) = \cap \{G \subseteq R(x) : G \cap G(f) \neq \varphi\}.$$

A function $f: U_1 \longrightarrow U_2$ is said to be rough in the approximation space $A = (U, R)$, where $A_1 = (U_1, R_1)$ and $A_2 = (U_2, R_2)$ are approximation spaces and $A = A_1 \times A_2, U = U_1 \times U_2$ if $\underline{R}(G(f)) \neq \overline{R}(G(f))$, otherwise $f$ is an exact function.

**Example 7.1** Let $U_1 = \{a,b,c,d,e\}$ and $U_2 = \{1,2,3,4,5,6\}$ be two universes, we define the function $f: U_1 \longrightarrow U_2$, by its graph $G(f) = \{(a,1), (a,6), (b,6), (c,5), (c,6), (e,6)\}$. Consider the blocks of the binary relations $R_1$ and $R_2$ as follows:

$$R_1(x) = \{\{a,c\}, \{a,b\}, \{d,e\}\},$$
$$R_2(x) = \{\{1,3\}, \{2,4,5\}, \{3,4\}, \{6\}\}.$$

Then

$$R(x) = R_1(x) \times R_2(x)$$
$$= \{\{(a,1),(a,3),(c,1),(c,3)\}, \{(a,2),(a,4),(a,5),$$
$$(c,2),(c,4),(c,5)\}, \{(a,3),(a,4),(c,3),(c,4)\}, \{(a,6),(c,6)\},$$
$$\{(a,1),(a,3),(b,1),(b,3)\}, \{(a,2),(a,4),(a,5),$$
$$(b,2),(b,4),(b,5)\}, \{(a,3),(a,4),(b,3),(b,4)\},$$
$$\{(a,6),(b,6)\}, \{(d,1),(d,3),(e,1),(e,3)\}, \{(d,2),(d,4),(d,5),$$
$$(e,2),(e,4),(e,5)\}, \{(d,3),(d,4),(e,3),(e,4)\}, \{(d,6),(e,6)\}\}$$

Then we have:
$$\underline{R}(G(f)) = \{(a,6), (b,6), (c,6)\},$$

$$\overline{R}(G(f)) = \{(a,1),(a,3),(c,1),(c,3),(a,6),(b,6),(c,6)$$
$$,(a,2),(a,4),(a,5),(c,2),(c,4),(c,5),(d,6),(e,6)\}$$

Therefore the function $f$ is a rough function such that $\underline{R}(G(f)) \neq \overline{R}(G(f))$.

When we have two approximation spaces defined by two equivalence relations, we have the following proposition that government the product space.

**Proposition 7.1** Let $A_1 = (U_1, R_1)$ and $A_2 = (U_2, R_2)$ be two arbitrary approximation spaces. Then we have: $(U_1 \times U_2)/R_1 \times R_2 = (U_1/R_1) \times (U_2/R_2)$.

**Proof**: suppose that: $u_1, u_2 \in U_1$ and $v_1, v_2 \in U_2$, then we have:
$$((u_1, v_1), (u_2, v_2)) \in R_1 \times R_2 \text{ iff } (u_1, u_2) \in R_1 \text{ and } (v_1, v_2) \in R_2.$$
Suppose again that: $[(u_1, v_1)]_{R_1 \times R2} \in (U_1 \times U_2)/R_1 \times R_2$.

Then we have 
$$[(u_1, v_1)]_{R_1 \times R2} = \{(u_2, v_2) : ((u_1, v_1), (u_2, v_2)) \in R_1 \times R_2\}$$
$$= \{(u_2, v_2) : (u_1, u_2) \in R_1, (v_1, v_2) \in R_2\},$$
$$= \{u_2 : (u_1, u_2) \in R_1\} \times \{v_2 : (v_1, v_2) \in R_2\}$$
$$= [u_1]_{R_1} \times [v_1]_{R_2}.$$

Then we have the result as: $(U_1 \times U_2)/R_1 \times R_2 = (U_1/R_1) \times (U_2/R_2)$. □

Let $f : (U_1, R_1) \to (U_2, R_2)$ be any function, where $A_1 = (U_1, R_1)$ and $A_2 = (U_2, R_2)$ are arbitrary approximation spaces. We define the relation $G(f) = \{(x, f(x)) : x \in U_1\}$ to be the graph of the function $f$. The rough approximations of $G(f)$ is defined as follows:

$$\underline{R}(G(f)) = \{(u_1, u_2) \in U_1 \times U_2 : [(u_1, u_2)]_R \subseteq G(f), R = R_1 \times R_2\},$$
$$\overline{R}(G(f)) = \{(u_1, u_2) \in U_1 \times U_2 : [(u_1, u_2)]_R \cap G(f) \neq \varphi, R = R_1 \times R_2\}.$$

Accordingly, the function $f$ is rough if $\underline{R}(G(f)) \neq \overline{R}(G(f))$, otherwise $f$ is an exact function. The pair $(\underline{R}(G(f)), \overline{R}(G(f)))$ is called a rough pair of relations.

The following theorems give the conditions on approximation spaces that give exact functions, one-to-one, surjective and continuous functions.

**Theorem 7.1** The function $f : U_1 \to U_2$ is an exact function for any selective approximation spaces $A_1 = (U_1, R_1)$ and $A_2 = (U_2, R_2)$.

**Proof:** The selective approximation space property means that, $[(u, v)]_R = \{(u, v)\}, u \in U_1, v \in U_2, R = R_1 \times R_2$. Then we have: $\underline{R}(G(f)) = \overline{R}(G(f))$, which yield to that the function $f$ is an exact function.□

**Theorem 7.2** The function $f : U_1 \to U_2$ is one-to-one function for any selective approximation spaces $A_1 = (U_1, R_1)$ and $A_2 = (U_2, R_2)$ if and only if both $\underline{R}(G(f))$ and $\overline{R}(G(f))$ are one-to-one functions.

**Proof:** The proof is directly using the definition of selective approximation space and using the technology in Theorem 3.1. □

**Theorem 7.3** The function $f: U_1 \to U_2$ is surjective function for any selective approximation spaces $A_1 = (U_1, R_1)$ and $A_2 = (U_2, R_2)$ if and only if both $\underline{R}(G(f))$ and $\overline{R}(G(f))$ are surjective functions.

**Proof:** As in the technique used in Theorem 3.1. □

**Theorem 7.4** The function $f: U_1 \to U_2$ is continuous function for any selective approximation spaces $A_1 = (U_1, R_1)$ and $A_2 = (U_2, R_2)$ if and only if both $\underline{R}(G(f))$ and $\overline{R}(G(f))$ are continuous functions.

**Proof:** As in the technique used in Theorem 7.2. □

When we have two topological spaces, generated using two bases $\beta_{R_1}, \beta_{R_2}$ where $A_1 = (U_1, R_1)$ and $A_2 = (U_2, R_2)$ are two approximation spaces. Then we have the following proposition that government the product topology.

**Proposition 7.2** Let $T_1 = (U_1, \tau_1)$ and $T_2 = (U_2, \tau_2)$ be two arbitrary topological spaces. Then we have: $(U_1 \times U_2) / \beta_{R_1} \times \beta_{R_2} = (U_1 / \beta_{R_1}) \times (U_2 / \beta_{R_2})$.

**Proof:** Similar as the proof of Proposition 7.1. □

The rough pairs of relations satisfied the following two important theorems.

**Theorem 7.5** For the quasi-discrete product topological space $(U_1 \times U_2, \tau)$, if $(\underline{R}(G(f)), \overline{R}(G(f)))$ is a rough pair of relations, and $(Q, \tau')$ is a subspace of $(U_1 \times U_2, \tau)$ such that $Q$ is closed in $\tau$, then $(\underline{R}(G(f)) \cap Q, \overline{R}(G(f)) \cap Q)$ is a relative rough pair of relations when $\underline{R}(G(f)), \overline{R}(G(f)), Q \subset U_1 \times U_2, \underline{R}(G(f)) \subset Q$.

**Proof:** The pair $(\underline{R}(G(f)), \overline{R}(G(f)))$ is a rough pair of relations in $(U_1 \times U_2, \tau)$, if the following condition satisfied:
1- $\underline{R}(G(f))$ is an open relation in $(U_1 \times U_2, \tau)$.
2- $\overline{R}(G(f))$ is a closed relation in $(U_1 \times U_2, \tau)$.
3- $\underline{R}(G(f)) \subset \overline{R}(G(f))$
4- The relation $\overline{\underline{R}(G(f)) - \overline{R}(G(f))}_\tau$ contains a relation $S$ of $U_1 \times U_2$ such that $S_\tau^o = \varphi$ and $\overline{R}(G(f)) - \overline{(\underline{R}(G(f))}_\tau \subset S_\tau$

Only we need to prove that $(\underline{R}(G(f)) \cap Q, \overline{R}(G(f)) \cap Q)$ is a rough pair of relations in $(Q, \tau')$. The proof will end by:

1- Since $\underline{R}(G(f))$ is an open relation in $(U_1 \times U_2, \tau)$, and $(Q, \tau')$ is a subspace of $(U_1 \times U_2, \tau)$, then $\underline{R}(G(f)) \cap Q$ is an open relation in $(Q, \tau')$.

2- Since $\overline{R}(G(f))$ is a closed relation in $(U_1 \times U_2, \tau)$, then there is an open relation $S$, such that $\overline{R}(G(f)) = U_1 \times U_2 - S$, then $\overline{R}(G(f)) \cap Q = (U_1 \times U_2) \cap Q - S \cap Q = Q - S \cap Q$, but $S \cap Q$ is an open relation in $(Q, \tau')$, then $\overline{R}(G(f)) \cap Q$ is a closed relation in the subspace $(Q, \tau')$.

3- Since $\underline{R}(G(f)) \subset \overline{R}(G(f))$, then $\underline{R}(G(f)) \cap Q \subset \overline{R}(G(f)) \cap Q$.

4- By selecting $S = R \cap Q, R = R_1 \times R_2$, then the relation $\overline{R(G(f))} \cap Q - \overline{(R(G(f)))}_\tau$ contains the relation $S$, and we need to prove the two sub conditions:
   a) $S_\tau^o = \varphi$,
   b) $\overline{R(G(f))} - \overline{(R(G(f)))}_\tau \subset \overline{S}_\tau$.

For the proof of Part a) $S_\tau^o = \varphi$, suppose that $S_\tau^o \neq \varphi$, then there is an $\tau$-open relation $G \subset Q$ such that $G \subset S$ but $S = R \cap Q$, i.e., $G \subset R$, but $G = G' \cap Q$ such that $G'$ is an open relation in $(U_1 \times U_2, \tau)$, then $G' \cap Q \subset R$ hence, $(G' \cap Q)_\tau^o \subset R_\tau^o$, but $R_\tau^o = \varphi$, that give contradiction, then must $S_\tau^o = \varphi$.

For the proof of Part b) $\overline{R(G(f))} - \overline{(R(G(f)))}_\tau \subset \overline{S}_\tau$.
Since $(\underline{R(G(f))}, \overline{R(G(f))})$ is a rough pair in $(U_1 \times U_2, \tau)$, then there is a relation $R' \subset U_1 \times U_2$, such that $\underline{R(G(f))} = R'^o_\tau$ and $\overline{R(G(f))} = \overline{R'}_\tau$ since, $R \subset \overline{R(G(f))} - \overline{(R(G(f)))}_\tau$, we have $S = R \cap Q = R' \cap Q - \overline{(R(G(f)))}_\tau$.
Now, let $(u,v) \in \overline{R(G(f))} \cap Q - \overline{(R(G(f)))}_\tau$, then $(u,v) \in \overline{R(G(f))} \cap Q$ and $(u,v) \notin \overline{(R(G(f)))}_\tau$.
Now if $(u,v) \in \overline{R(G(f))} \cap Q$, then $(u,v) \in S$, and $(u,v) \in \overline{S}_\tau$.

Finally, if $(u,v) \notin R' \cap Q$ and $(u,v) \in \overline{R(G(f))} \cap Q$ and $(u,v) \notin \overline{(R(G(f)))}_\tau$, hence $(u,v) \in \overline{R(G(f))}$ and $(u,v) \in Q$. Now $(u,v) \in \overline{R'}_\tau$, then there is an open relation $G$ in $\tau$ such that $(u,v) \in G$ and $G \cap R' \neq \varphi$, but $(u,v) \notin \overline{(R(G(f)))}_\tau$, then $(u,v) \in G - \overline{(R(G(f)))}_\tau = G \cap [\overline{(R(G(f)))}_\tau]^c$. But since $\underline{R(G(f))} = R'^o_\tau$ is an open relation in $\tau$, and $\underline{R(G(f))} = \underline{R(G(f))} \cap Q$ is an open relation in $\tau'$, then $\overline{(R(G(f)))}_\tau$ is a closed relation in $\tau$, hence $[\overline{(R(G(f)))}_\tau]^c$ is an open relation in $\tau$, hence $G \cap [\overline{(R(G(f)))}_\tau]^c$ is an open relation containing $(u,v)$, then $G \cap [\overline{(R(G(f)))}_\tau]^c \cap R' \neq \varphi$. This yields to $G \cap [R' - \overline{(R(G(f)))}_\tau] \neq \varphi$, i.e. $G \cap [R' - \overline{(R(G(f)))}_\tau] \cap Q \neq \varphi$, such that $(u,v) \in Q$. Hence $G \cap [R' \cap Q - \overline{(R(G(f)))}_\tau \cap Q] \neq \varphi$, but $\overline{R(G(f))} \subset Q$, then $G \cap [R' \cap Q - \overline{(R(G(f)))}_\tau] \neq \varphi$. But we have $R = R' - \overline{(R(G(f)))}_\tau$, hence $R \cap Q = R' \cap Q - \overline{(R(G(f)))}_\tau$ implies that $S = R' \cap Q - \overline{(R(G(f)))}_\tau$. Then $G \cap S \neq \varphi$, but $(u,v) \in Q$, hence $(G \cap Q) \cap S \neq \varphi$. But $G \cap Q = G'$ and $G'$ is an open relation in $\tau'$ that containing $(u,v)$, then $G \cap S \neq \varphi$, hence $(u,v) \in \overline{S}_\tau$. □

**Theorem 7.6** Let $(\underline{R(G(f))}, \overline{R(G(f))})$ be a rough pair of relations in the product topological space $(U_1 \times U_2, \tau)$, and let $(Q, \tau')$ be a subspace of $(U_1 \times U_2, \tau)$ such that $Q$ is any relation of $U_1 \times U_2$. Then there is a relation $P \subset Q$ such that $\underline{R(G(f))} \cap Q = P_\tau^o$ and $\overline{R(G(f))} \cap Q = \overline{P}_\tau$.

**Proof:** We can define $P = \underline{R(G(f))} \cap Q$, then $P \subset Q$, then $(\underline{R(G(f))} \cap Q)_\tau^o = P_\tau^o$, but $\underline{R(G(f))} \cap Q$ is an open relation in $(Q, \tau')$ i.e., $(\underline{R(G(f))} \cap Q)_\tau^o = \underline{R(G(f))} \cap Q$, hence $\underline{R(G(f))} \cap Q = P_\tau^o$.

Finally for $\overline{R(G(f))} \cap Q = \overline{P}_\tau$, since $(\underline{R(G(f))}, \overline{R(G(f))})$ is a rough pair in $(U_1 \times U_2, \tau)$, then $\underline{R(G(f))} \cap Q \subset \overline{R(G(f))} \cap Q$ i.e., $P \subset \overline{R(G(f))} \cap Q$ implies that $\overline{P}_\tau \subset \overline{\overline{R(G(f))} \cap Q}_\tau$. But $\overline{R(G(f))} \cap Q$ is a closed relation, i.e., $\overline{\overline{R(G(f))} \cap Q}_\tau = \overline{R(G(f))} \cap Q$, hence $\overline{P}_\tau \subset \overline{R(G(f))} \cap Q$.

For $\overline{R}(G(f)) \cap Q \subset \overline{P}_\tau$, let $(u,v) \in \overline{R}(G(f)) \cap Q$, then $(u,v) \in \overline{R}(G(f))$ and $(u,v) \in Q$, but since $\underline{R}(G(f)) \cap Q = P$, then $Q \subset [\underline{R}(G(f)) - P]^c = U_1 \times U_2 - [\underline{R}(G(f)) - P]$, hence $(u,v) \in U_1 \times U_2 - [\underline{R}(G(f)) - P]$ i.e., $(u,v) \in P$, hence $(u,v) \in \overline{P}_\tau$, then $\overline{R}(G(f)) \cap Q \subset \overline{P}_\tau$. Then we have $\overline{R}(G(f)) \cap Q = \overline{P}_\tau$. □

Let $(U_1 \times U_2, \tau)$ be a product space. For any relation $Q \subset U_1 \times U_2$ define the subspace $(Q, \tau')$ of $(U_1 \times U_2, \tau)$. We define the equivalence relation $E(\tau')$ on the power set $P(Q)$ by: $(R_1, R_2) \in E(\tau') \Leftrightarrow (R_1)^o_{\tau'} = (R_2)^o_{\tau'}, (\overline{R_1})_{\tau'} = (\overline{R_2})_{\tau'}$ for any $R_1, R_2 \in P(Q)$. The set $P(Q)/E(\tau')$ is a partition of $P(Q)$ and any class $\eta \in P(Q)/E(\tau')$ is called a relative topological rough relations.

**Theorem 7.7** For any product topological space $(U_1 \times U_2, \tau)$ and for any subspace $(Q, \tau')$ of it the function $f : P(Q)/E(\tau') \to \eta'$, defined by: $f(R) = ((R)^o_{\tau'}, (\overline{R})_{\tau'}), R \in \eta'$ is bijection, where $\eta'$ is the set of all relative rough pairs.

**Proof**: the proof is directly by Theorems 7.2 & 7.3. □

Let $(U_1, \tau_1)$ and $(U_2, \tau_2)$ be any two topological spaces, where $\beta_1$ and $\beta_2$ are any two bases for $\tau_1$ and $\tau_2$. Then we define the base $\beta = \beta_1 \times \beta_2$ of the topology $\tau = \tau_1 \times \tau_2$.

We define the approximations for any subset $H \subseteq U_1 \times U_2$:

$$(H)^o_\tau = \cup \{G \subseteq \beta : G \in H\},$$
$$\overline{(H)}_\tau = \cap \{G \subseteq \beta : G \cap H \neq \varphi\}.$$

The function $f$ on $U_1 \times U_2$ is called a topological rough continuous function at the point $(x,y) \in U = U_1 \times U_2$ if $f^{-1}(V(f(x,y))) \subseteq \tau$ for all open set $V(f(x,y)) \in \tau$. The function $f$ is topological rough continuous on $U_1 \times U_2$ if it is a topological roughly continuous at every point of $U_1 \times U_2$.

**Example 7.2** Consider the topology $\tau_1 = \{U_1, \varphi, \{a\}, \{b,c,d\}\}$, on $U_1 = \{a,b,c\}$ and the topology $\tau_2 = \{U_2, \varphi, \{3\}, \{1,2,4\}\}$ on $U_2 = \{1,2,3,4\}$. The bases $\beta_1 = \{\{a\}, \{b,c,d\}\}$ and $\beta_2 = \{\{3\}, \{1,2,4\}\}$ are of $\tau_1$ and $\tau_2$ respectively.

We defined the function $f : U_1 \times U_2 \to U_1 \times U_2$ as follows:
$$f(x,y) = (a,3),$$

Then we have:
$\beta = \beta_1 \times \beta_2 = \{\{a\}, \{b,c,d\}\} \times \{\{3\}, \{1,2,4\}\}$
$= \{\{a\} \times \{3\}, \{a\} \times \{1,2,4\}, \{b,c,d\} \times \{3\}, \{b,c,d\} \times \{1,2,4\}\}$
$\tau = \tau_1 \times \tau_2 = \{U_1, \varphi, \{a\}, \{b,c,d\}\} \times \{U_2, \varphi, \{3\}, \{1,2,4\}\}$
$= \{U_1 \times U_2, U_1 \times \varphi, U_1 \times \{3\}, U_1 \times \{1,2,4\},$
$\varphi \times U_2, \varphi \times \varphi, \varphi \times \{3\}, \varphi \times \{1,2,4\},$
$\{a\} \times U_2, \{a\} \times \varphi, \{a\} \times \{3\}, \{a\} \times \{1,2,4\},$
$\{b,c,d\} \times U_2, \{b,c,d\} \times \varphi, \{b,c,d\} \times \{3\}, \{b,c,d\} \times \{1,2,4\}\}$

Then for any point $(x,y) \in U_1 \times U_2$ we have $f(x,y) = (a,3)$, then all open sets containing $(a,3)$ are:
$V_1 = \{a\} \times U_2 = \{(a,1), (a,2), (a,3), (a,4)\}$,
$V_2 = \{a\} \times \{3\} = \{(a,3)\}$
Then the inverse function of these open sets is:

$f^{-1}(V_1) = U_1 \times U_2$,
$f^{-1}(V_2) = U_1 \times U_2$

Then the function $f$ is topological rough continuous at every point of $U_1 \times U_2$.

## 8. Future Applications of Topological Rough Functions on Information Systems

In this section, we will define a function between two information systems and give all needed conditions for them. Functions on information system can produce the reducts and the core of this system by the projection the system on subsystems. We will define the image of rough set using some types of these functions. Finally, we define the topological rough functions on information systems and study some of their properties.

The reader can review about information systems in [15, 18] to know about the structure and the types and the different methods of reduction of information systems.

### 8.1 Suggested Functions on an Information System

Suppose an information system $T = (U, C, D)$ where $U$ is the set of objects of this system (Patients, Plants,…). $C$ is the condition attributes of these objects (Temperature, Muscle pain,…). $D$ is the expert decisions about the condition attribute that objects suffer from.

We define the projection (restriction) function $f_c : P(C) \times P(C) \to P(C)$, where $P(C)$ is the power set of the condition attributes as follows:

$$f_c((B, B')) = \begin{cases} C & \text{,if } POS_B(D) \neq POS_{B'}(D), \forall B' \subseteq C \\ B' & \text{, if } POS_B(D) = POS_{B'}(D), \forall B' \subset B \end{cases}$$

Figure 1 below gives an example for a projection function on information system. The core of such systems is given by taking the intersection of all these projection functions on that system.

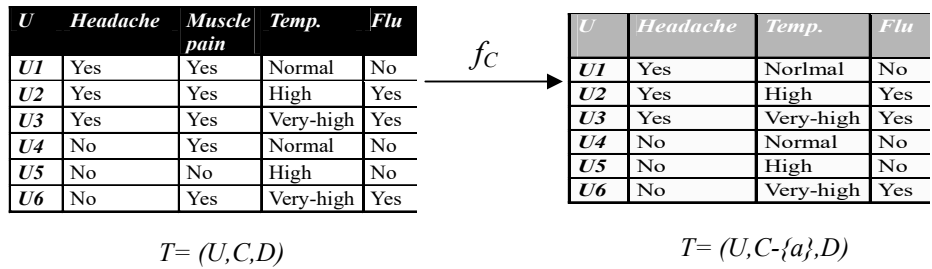

T= (U,C,D)    T= (U,C-{a},D)

**Figure 1**: Some reducts of information system by projection function

The topological rough, continuous functions of information systems can be defined as follows:

The function $f : (U, C, D) \to (U, C', D)$ is a topological roughly continuous of the object $x \in U$ if $f^{-1}(D_{(C-a)}(x)) = D_C(x)$, where $D_C(x) = \{y \in U : D_C(x) = D_C(y)\}$. The function $f$ is topological roughly continuous on $U$ if it is topological rough continuous for every object of $U$.

### 8.2 Suggested Functions on orders information systems

By a discernibility matrix of information system $T$, denoted $M(T)$, we will mean $n \times n$ matrix defined as follows:
$M(T) = \{m_{ij} : i, j = 1, 2, 3, ..., n\}$ where,

$$m_{ij} = \begin{cases} \{a \in C : a(x_i) \neq a(x_j)\} : if\ \exists b \in D, b(x_i) \neq b(x_j) \\ \lambda \qquad\qquad\qquad\qquad\quad : if\ \forall b \in D, b(x_i) \neq b(x_j) \end{cases}$$

Such that $a(x_i)$ or $a(x_j)$ belongs to the C-positive region of $D$. $m_{ij}$ is the set of all condition attribute that classify objects $a(x_i)$ and $a(x_j)$ into different classes. $m_{ij} = \lambda$ denotes that this case does not need to be considered.

The discernibility function $f : T = (U, C, D) \to M(T)$ of an information system is defined as follows:

For any object $x_i \in U$: $f_T(x_i) = \wedge_j \{\vee m_{ij} : i \neq j, j \in \{i, 2, ..., n\}\}$, where $\vee m_{ij}$ is the disjunction of all variables $b \in m_{ij}$ when $m_{ij} \neq \varphi$ and $\vee m_{ij} = 0$ when $m_{ij} = \varphi$ and $\vee m_{ij} = 1$ when $m_{ij} = \lambda$.

Figure 2 below gives an example for a discernibility function on information system.

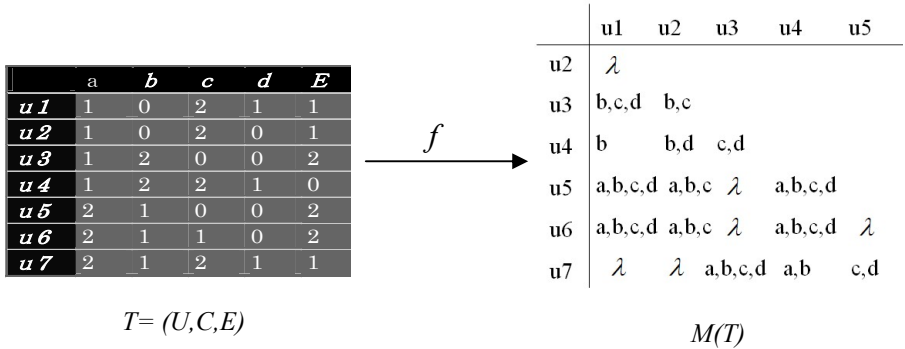

Figure 2: Discernibility function of information system

According to Figure 2 the function $f$ transfer the system $T = (U, C, E)$ into the Discernibility $M(T)$ and the reduction of this system can be obtained as follows:
$$f_T(x_i) = f_T(a, b, c, d) = b \wedge (a \vee b) \wedge (c \vee d) \wedge (b \vee d) \wedge$$
$$(a \vee b \vee c) \wedge (a \vee b \vee c \vee d)$$

Then we have:
$$f_T(x_i) = b \wedge (c \vee d)$$

Accordingly, the system $T = (U, C, E)$ has two reducts namely $R_1 = \{b, c\}$ and $R_2 = \{b, d\}$ with core $CORE(T) = \{b\}$.

## 9. Predictions of Patients Classification Data Using Rough Continuous Functions

Our aim in this application problem to find the future recommendations for patients show them appropriately greeted approach combines treatment and exercise to reach results explain the function of each presentation of the positive and negative impact on the patient.

The decision of the Physician, according to the medical reports is the continuation of the medical tests are all for another or off the medical analyzes the patient's condition is stable loft insensitively to healthy style workout constantly.

The structure $S = (U, At, \{V_a : a \in At\}, f_a, \{R_P : P \subseteq At\})$ is the mathematical style of information system of our patients problem. The set $U$ is the system universe that we selected to be a set of five patients. The set $At$ is the set of attributes of these patient's namely tests functions such as liver, kidney, and heart functions and so on. The set $V_a$ is values of each attribute $a \in At$. Finally, $f_a : U \to V_a$ is the information function such that $f_a(x) \in V_a$.

For any subset $B \in At$, we define the relation:

$R_P = \{(x,y): |f_a(x) - f_a(y)| < \alpha, a \in P, x, y \in U, \alpha \in R^+\}$, for $a \in At$, we define the class $A_{R_a}$ as follows: $A_{R_a} = \{R_a(x): x \in U\}$, where $R_a(x) = \{y: xR_ay\}$.

The structure $DS = (U, \{At, D\}, \{V_a : a \in At\}, f_a, \{R_P : P \subseteq At\})$ is a decision table, where $D$ is the set of decisions, that represents for each patient if he need surgery or enough need drugs.

We define the relation of the decision attribute $D$ by:

$R_D = \{(x,y): f_a(x) = f_a(y), a \in D, x, y \in U\}$.

The class of this relation is $R_D(x) = \{y: xR_Dy\}$. The set of all classes is $A_{R_a} = \{R_a(x): x \in U\}$. We define the set $P \subseteq At$ to be a reduct of $At$, if $\tau_D \subseteq \tau_P$ and $P$ is a minimal.

Basic data of five patients before the surgery are given in Table 2 (The decision System of Patients). Each patient will measure these medical functions periodical every three months. After a period of time we need to predict the results of the medical tests of patients in any time and accordingly they can stop drugs. Therefore, we defined the prediction function $f_P : DS \to \overline{\overline{DS}}$, where $\overline{\overline{DS}}$ is the decision system of patients after time $t$ (Dynamic Decision System of Patients).

| Patients ID | Liver Functions | | Kidney Functions | Heart Efficiency | Decision |
|---|---|---|---|---|---|
| | A1 | A2 | A3 | A4 | D |
| X1 | 35 | 45 | 6.8 | 412 | Need Surgery |
| X2 | 45 | 44 | 4.2 | 420 | Need Drugs |
| X3 | 42 | 38 | 5.8 | 430 | Need Surgery |
| X4 | 30 | 44 | 9.7 | 480 | Need Surgery |
| X5 | 36 | 32 | 5.4 | 450 | Need Drugs |

Table 2: The infection information system of patients

Now if we choose for the liver function attributes $P_1 = LF = \{A_1, A_2\}$ the threshold $\alpha_1 = 4$, then $R_{P_1}(U) = \{\{X3\}, \{X5\}, \{X1, X4, X5\}, \{X3, X5\}, \{X2, X3, X4\}\}$. The topology generated by $R_{P_1}$ is given by:

$\tau_{P_1} = \{U, \varphi, \{X3\}, \{X5\}, \{X1, X4, X5\}, \{X3, X5\}, \{X2, X3, X4\},$
$\{X1, X3, X4, X5\}, \{X2, X3, X4, X5\}\}$

For kidney functions we can choose $\alpha_1 = 2.5$ for $P_2 = KF = \{A_3\}$, then
$R_{P_2}(U) = \{\{X4\}, \{X1, X4\}, \{X1, X2, X3, X5\}\}$,
$\tau_{P_2} = \{U, \varphi, \{X4\}, \{X1, X4\}, \{X1, X2, X3, X5\}\}$.

For the heart efficiency attribute $P_3 = HE = \{A_4\}$ we can choose $\alpha_3 = 20$, then
$R_{P_3}(U) = \{\{X4, X5\}, \{X1, X2, X3, X5\}, \{X1, X2, X3, X4\}\}$,
$\tau_{P_3} = \{\{X4, X5\}, \{X1, X2, X3\}, \{X1, X2, X3, X5\},$
$\{X1, X2, X3, X4\}\}$

For the decision attribute we have $R_D(U) = \{\{X1, X3, X4\}, \{X2, X5\}\}$, then the topology of decisions is $\tau_D = \{U, \varphi, \{X1, X3, X4\}, \{X2, X5\}\}$.

The condition attributes are exactly three attributes namely $At = \{LF, KF, HE\}$. The numbers of non-trivial subsets of the set of all condition attributes are seven subsets namely $\{P1, P2, P3, \{P1, P2\}, \{P1, P3\}, \{P2, P3\}, \{P1, P2, P3\}\}$.

Now we will calculate the classes of the residue subsets by taking the intersections as follows:

$R_{P_1,P_2}(U) = R_{P_1}(U) \cap R_{P_2}(U) = \varphi$, with topology: $\tau_{P_1,P_2} = \{U, \varphi\}$

$R_{P_1,P_3}(U) = R_{P_1}(U) \cap R_{P_3}(U) = \varphi$, with topology: $\tau_{P_1,P_3} = \{U, \varphi\}$

$R_{P_2,P_3}(U) = R_{P_2}(U) \cap R_{P_3}(U) = \{\{X1, X2, X3, X5\}\}$, with topology:

$\tau_{P_2,P_3} = \{U, \varphi, \{X1, X2, X3, X5\}\}$

The covering class of universe using all condition attributes is giving as follows:

$R_{P_1,P_2,P_3}(U) = R_{P_1}(U) \cap R_{P_2}(U) \cap R_{P_3}(U) = \varphi$, with topology: $\tau_{P_1,P_3} = \{U, \varphi\}$. Then the system given in Table 2 has no topological reducts.

Now we define the function $f_P : DS \to \overline{\overline{DS}}, P \subseteq At$ by $f_P(Xi) = Xi, i = 1, 2, 3, 4, 5$. Then according the function $f_P(Xi)$ the image of Table 2 after a period of three months has no topological reducts and this function is one to one rough continuous function.

## 10. Conclusion and Future Work

The emergence of topology in the construction of some rough functions will be the bridge for many applications and will discover the hidden relations between data. Topological generalizations for the concept of rough functions opens the way for connecting rough continuity with the area of near continuous functions.

Applications of topological rough functions on information systems open the door about many transformations among different types of information systems such as multi-valued and single valued information systems.

Future applications of our approach in computer can be as follows:
- In information Retrieval Fields: we can modify the query running online by define a function that coverted documents to weighted vectors of words of that document. Then we can extract the results of weights in a decision table that we can classify the documents accordingly the reduction of this table. The query is constructed by defining a Boolean function of all words of the reduction.
- Classification and Summarization of Documentation using topological functions of neighborhood systems defined on documents.